\RequirePackage{fix-cm}
\documentclass[smallcondensed]{svjour3}       
\smartqed  
\usepackage{epsfig,amssymb,amsmath}
\usepackage{amsfonts}
\usepackage{graphicx}

\newcommand{\sta}{{\rm sta}}
\newcommand{\ext}{{\rm ext}}
\newcommand{\Diag}{{\rm Diag}}

\newcommand{\eb}{\begin{equation}}
\newcommand{\ee}{\end{equation}}
\newcommand{\la}{\langle}
\newcommand{\ra}{\rangle}
\newcommand{\bx}{\boldsymbol{x}}
\newcommand{\by}{\boldsymbol{y}}
\newcommand{\bvsig}{\boldsymbol{\varsigma}}
\newcommand{\bveps}{\boldsymbol{\varepsilon}}
\newcommand{\barbx}{\boldsymbol{\bar{x}}}

\newcommand{\barby}{\boldsymbol{\bar{y}}}
\newcommand{\barvsig}{\bar{\bvsig}}
\newcommand{\PP}{{\Pi}}
\newcommand{\PPd}{\PP^{d}}
\newcommand{\half}{\frac{1}{2}}
\newcommand{\bff}{\boldsymbol{f}}
\newcommand{\Rn}{\mathbb{R}^{n}}
\newcommand{\Rm}{\mathbb{R}^m}
\newcommand{\real}{\mathbb{R}}
\newcommand{\calP}{\mathcal{P}}

\newcommand{\R}{\mathbb{R}}
\newcommand{\Lbd}{\Lambda}
\newcommand{\Lam}{\Lambda}
\newcommand{\barLam}{\bar{\Lambda}}
\newcommand{\bxi}{\boldsymbol{\xi}}
\newcommand{\calV}{\mathcal{V}}
\newcommand{\calS}{\mathcal{S}}
\newcommand{\vsig}{\varsigma}
\newcommand{\bb}{\boldsymbol{b}}
\newcommand{\bp}{\boldsymbol{p}}

\newcommand{\bq}{\boldsymbol{q}}
\newcommand{\Gap}{G_{ap}}
\newcommand{\Pf}{P_\flat}
\newcommand{\Ps}{P_\sharp}
\newcommand{\Qf}{Q_\flat}
\newcommand{\Qs}{Q_\sharp}
\newcommand{\barbvsig}{\bar{\boldsymbol{\varsigma}}}
\newcommand{\calX}{\mathcal{X}}
\newcommand{\calQ}{\mathcal{Q}}
\newcommand{\calY}{\mathcal{Y}}
\newcommand{\calXs}{\calX_\sharp}
\newcommand{\calSs}{\calS_\sharp}
\newcommand{\VV}{V}
\newcommand{\UU}{U}
\newcommand\barW{{\bar W}}
\newcommand\barF{{\bar F}}

\newcommand{\BB}{B}
\newcommand{\TT}{T}
\newcommand{\DD}{D}
\newcommand{\eps}{\epsilon}
\newcommand{\calU}{{\cal U}}
\newcommand{\calE}{{\cal E}}
\newcommand{\uu}{u}
\newcommand{\baru}{{\bar u}}

\newcommand{\G}{G\^{a}teaux }
 \journalname{J Glob Optim}
\begin{document}

\title{Triality Theory for General  Unconstrained
Global Optimization Problems}


\author{David Y. Gao        \&  Changzhi Wu
}


\institute{D.Y. Gao \at
              School of Science, Information Technology and Engineering, University of Ballarat, Victoria 3353, Australia \\
              \email{d.gao@ballarat.edu.au}           
           \and
           C.Z. Wu \at
            School of Science, Information Technology and Engineering, University of Ballarat, Victoria 3353, Australia \\
              \email{c.wu@ballarat.edu.au}
}

\date{Received: date / Accepted: date}
\maketitle

\begin{abstract}
 Triality theory is proved for a  general unconstrained
global optimization problem. The method adopted is simple but mathematically rigorous.
 Results show that if the primal problem and its canonical dual have the same dimension,
the triality theory holds strongly in the tri-duality form as it was originally proposed.
Otherwise, both
the canonical min-max duality and the double-max duality still hold strongly, but
the double-min duality holds weakly in a super-symmetrical form as it was expected.
Additionally, a complementary weak saddle min-max duality theorem is discovered.
Therefore, an open problem on this statement left in 2003 is solved completely.
This theory can be used to identify not only the global minimum, but also
the largest local minimum,   maximum, and saddle points. Application is
 illustrated. Some fundamental concepts  in  optimization
and remaining challenging problems in canonical duality theory   are discussed. \vspace{.5cm}
\end{abstract}



\section{Introduction}
The general global optimization problem to be solved is proposed in the following form
 \eb
 (\calP): \;\;  \ext \left\{\PP(\bx)= W(  \bx) + \half \la \bx ,
 A\bx \ra - \la \bx , \bff \ra \;|\;\bx\in \Rn\right\}, \label{PP0}
 \ee
where
$W(\bx)$ is a nonconvex  function,  $A\in\R^{n\times n}$ is a given symmetric matrix,
$\bff\in\Rn$ is a given  vector (source),
$\la *  ,  * \ra $
is an inner product in  $\Rn$, and
the  notation $\ext \{ * \}$ stands for finding global extrema of the function given in $\{ * \}$,
including both global minimum and  the largest local minimum and  maximum.
In order to have this general problem making sense in reality, the
  nonconvex function $W(\bx)$ should obey  certain fundamental rules in systems theory.

 Objectivity is a basic concept in science, which  is often attributed with the property of scientific measurements that can be measured independently of the observer.
 General description of the objectivity can be easily found on internet and in many mathematical physics  textbooks  (see \cite{holz,ogden}).
 Mathematical definitions of the objective set and objective function are given in the book \cite{GaoBook} (Chapter 6, page 288).
 Let
 \[
 {\cal Q} = \{ Q \in \real^{m\times m} | \; Q^T = Q^{-1}, \;\; \det Q = 1 \}
 \]
 be a proper
 orthogonal rotation group.
\begin{definition}[Objectivity and Isotropy]
{  A subset $\calY_a \subset \real^m$ is said to be {\em objective} if
$ Q \by \in \calY_a \;\; \forall \by \in \calY_a   \;  \mbox{ and } \;  \forall Q \in {\cal Q}.
$   A real-valued function $\TT:\calY_a  \rightarrow \real$ is said to be
 {\em objective}  if its domain is objective and
 \eb
\TT( Q\by) = \TT(\by) \;\; \forall \by \in \calY_a  \; \mbox{ and } \; \forall Q \in {\cal Q}.
 \ee

  A subset $\calY_a \subset \real^m$ is said to be {\em isotropic} if
$  \by Q^T \in \calY_a  \;\; \forall \by \in \calY_a  \; \mbox{ and } \; \forall Q \in {\cal Q}.
$
  A real-valued function $\TT :\calY_a  \rightarrow \real$ is said to be
 {\em isotropic}  if its domain is isotropic  and
 \eb
\TT(  \by Q^T) = \TT (\by) \;\; \forall \by \in \calY_a  \; \mbox{ and } \;  \forall Q \in {\cal Q}.
 \ee
 }\end{definition}

Geometrically speaking, the objectivity  means that
the function $\TT (\by)$ does not depend on  rotation, but
on certain measure (norm) of its variable $\by$.
Therefore, the most simple objective function is the $l_2$-norm
$\TT (\by) = \|\by \|$ since
$\|Q \by \|^2 =  \by^T Q^T Q \by = \by^T \by = \| \by\|^2 \;\; \forall Q \in \calQ$.
While the isotropy implies  that the function $\TT(\by)$ possesses a certain symmetry.
By the fact that
     $( \bx Q^T) (\bx Q^T)^T  = \bx \bx^T   \succeq 0 \;\; \forall  Q \in {\cal Q}$,
    the concept of isotropy plays important role in Semi-Definite Programming (SDP) and
    integer programming  \cite{gao-jimo07,gao-ruan-jogo10}.

The objectivity in science is also refereed as {\em frame invariance},   which
lays a foundation for mathematical physics and
systems theory.
In fact, the canonical duality theory was originally developed from this  concept
\cite{GaoBook}, which is the  reason  why this theory can be applied not only for modeling and
  analysis of complex systems, but also for solving a large class of nonconvex/nonsmooth/discrete
  problems in both mathematical physics and global optimization.
   In this
paper, we shall need only the following weak assumptions for the nonconvex function $W(\bx)$.
\begin{verse}
    \item[(A1).] The nonconvex function $W(\bx)$ is twice continuously differentiable.\\
    \item[(A2).] There exits a {\em geometrical operator}
  \eb
  \Lbd(\bx) = \left\{ \half \bx^T \BB^k \bx +\bb_k^T\bx
\right\} :\;\Rn \rightarrow \Rm
\ee
and a strictly convex
function $V:\real^m \rightarrow \real$ such that
 \eb
 W(\bx) = V(\Lbd (\bx)), \label{canonical}
 \ee
 where $\BB^k\in\mathbb{R}^{n\times n}$ and $\bb_k\in\Rn,
 k=1,\cdots,m$.\\
    \item[(A3).] The critical points of problem ($\calP$) are
    non-singular, i.e, if $\nabla\PP(\barbx)=0$, then $\det(\nabla^2\PP(\barbx))\neq
    0$.
\end{verse}

Based on Assumption (A2), the general problem (\ref{PP0}) can be reformulated in the following
canonical form:
\eb
(\calP): \;\;
\ext \left\{\PP(\bx)= \VV(\Lam( \bx)) + \half \la \bx ,
 A\bx \ra - \la \bx , \bff \ra \;|\;\bx\in \Rn\right\}. \label{PP}
 \ee
This  problem  arises extensively in many fields of
engineering and sciences, including  Euclidean distance geometry \cite{Floudas1,gao-ruan-pardalos},
 computational biology \cite{Floudas2,Ogden,zhang-gao}, numerical methods for solving a large
  class of nonconvex variational problems in mathematical physics
 \cite{Gao-amma03,Gao-Yu2008,santo-gao}, and much more.

Actually, the assumption (A2) is   the so-called {\em canonical transformation}
introduced in \cite{GaoBook}. The  idea of this transformation was from
Gao and Strang's original work \cite{GaoStrang89} on
nonconvex variational problems in large deformation theory, where
 the geometrical operator
   $\Lam({u}) = \half (\nabla {  u})^T (\nabla {  u})$ is a Cauchy-Riemann metric tensor field,
   which is an objective measure of the deformation gradient $\eps = \nabla u$,
   and $W(\nabla u) = V(\Lam(u))$ is a stored strain energy.
  By using finite element discretization  for the deformation  field $u(\bx)$, the
   nonconvex variational problems in infinite dimensional space can be reduced to the
   canonical global optimization problem $(\calP)$ (see \cite{Gao-Yu2008,santo-gao}).
       It is known in continuum physics that the stored energy $W$ is usually a nonconvex function of
   the linear measure $\nabla u$ (which is not a strain measure),
    but $V(e)$ is convex in term of the objective  measure $e = \Lam(u)$.
    Therefore, by this quadratic objective operator $\Lam(u)$, a
    {\em complementary gap function}
    was discovered by Gao and Strang in nonconvex variational analysis, and by which,
      complementary variational principles were recovered in fully nonlinear  equilibrium problems of
mathematical physics\footnote{In continuum physics, complementary variational principle means
perfect duality since any duality gap will violate certain physical laws. The existence of a complementary
variational principle was a well-known
debate existing for several decades in large  deformation theory  (see \cite{li-gupta}).
This problem was partially solved by Gao and Strang's work, and solved completely in 1999 \cite{gao-mecc}.}.
     They also proved that the nonnegative  gap function can be used to identify global minimizer of
     the nonconvex problem. Seven years later, it was discovered   that
     the negative gap function can be used to identify the largest local minimum and maximum.
     Therefore, the {\em triality theory} was first proposed in nonconvex mechanics \cite{gao-tri96,gao-amr}, and then generalized to global optimization \cite{GaoJOG00}.
     This triality theory is composed of a canonical  min-max duality and two pairs of double-min, double-max dualities, which reveals an intrinsic duality pattern in complex systems
     and has been used successfully for  solving a wide class of
     challenging problems in nonconvex analysis and global optimization \cite{GaoBook}.
     However, it was realized  in 2003 \cite{gao-opt03,Gao-amma03} that the double-min duality holds conditionally under  ``certain additional conditions".
 Recently, this problem is partly solved for a class of fourth order
polynomial optimization problems \cite{GaoWu,silva-gao-jmaa1}.

 The aim of this paper is to prove the
  triality theory  for the  general nonconvex  global optimization problem $(\calP)$.
   In the following sections, we first provide a brief review on the
   canonical duality theory and the associated triality theory.
   We will show that by the canonical transformation, the  nonconvex primal problem $(\calP)$
   can be reformulated as a canonical dual problem  without duality gap.
    Section 3 presents a  strong triality theory
  for the case that the primal problem and its canonical dual have the same dimension, i.e.
  $n=m$. We then show in Section 4 that this theory holds weakly
  for the case $n\neq m$. The ``certain additional conditions" for the double-min duality
  are provided.
  Application is illustrated  in Section 5.
  The original definition of Lagrangian, Lagrangian duality and its difference with the canonical duality  are discussed in Section 6.
  The paper is ended with some conclusion remarks and challenging  problems.

\section{Canonical Duality, Triality, and Open Problem}

  Let
\[
\calV_a = \{ \bxi \in \real^m | \;\; \bxi  = \Lam(\bx)\;\; \forall \bx \in \real^n \},
\]
\[
\calV^*_a = \{ \bvsig \in \real^m | \;\; \bvsig   = \nabla V(\bxi) \;\; \forall \bxi \in \calV_a \}.
\]
 By  (A1) and (A2) we know that $V:\calV_a \rightarrow \real$ is also a
twice continuously differentiable. Therefore,
its  Legendre conjugate $V^*:\calV^*_a\rightarrow \real $ can be uniquely defined
as
 \eb
 V^*(\bvsig) = \sta\left\{\la \bxi ; \bvsig\ra-V(\bxi)\;|\;\bxi\in\calV_a\right\},
 \ee
 where $\la * ;  * \ra $
is an inner product in  $\Rm$ and  $\sta\{\;\}$ stands for finding
stationary value of the expression given in  $\{\;\}$.
It is easy to verify that the   canonical duality relations
 \eb
\bvsig = \nabla V(\bxi) \;\;  \Leftrightarrow \;\; \bxi = \nabla
V^*(\bvsig)  \;\; \Leftrightarrow \;\; V(\bxi ) + V^*(\bvsig) =
\la \bxi ; \bvsig \ra
 \ee
  hold on $  \calV_a \times \calV^*_a$.

   Substituting
  $V(\Lbd(\bx))= \la \Lbd(\bx) ; \bvsig\ra - V^*(\bvsig)$, the
  primal function $\Pi(\bx)$ can be reformulated as  the total-complementary function
\cite{GaoBook}
 \eb
 \Xi(\bx,\bvsig) = \half \la  \bx ,  G(\bvsig) \bx \ra -
V^*(\bvsig) - \la \bx ,  F(\bvsig) \ra ,  \label{Complementary}
 \ee
where
\[
G(\bvsig) = A + \sum_{k=1}^m \vsig_k B^k , \;\;\;
 F(\bvsig)
= \bff - \sum_{k=1}^m\vsig_k \bb_k.
\]
For a fixed $\bvsig$, the criticality condition
$\nabla_{\bx}\Xi(\bx,\bvsig) = 0$ leads to the following
canonical equilibrium equation
 \eb
 G(\bvsig)\bx = F(\bvsig),  \label{Equi}
 \ee
which can be solved analytically to obtain\footnote{In
this paper $G^{-1}$ should be understood as a  generalized
inverse if $\det G = 0$ \cite{GaoJOG00}}
 $\bx = [G(\bvsig)]^{-1} F(\bvsig)$ for all $\bvsig$ in the canonical dual feasible space $\calS_a $
 defined by
\[
\calS_a = \{ \bvsig \in \calV^*_a | \; F(\bvsig)\in
\mathcal{C}_{ol}\left( G \left( \bvsig\right) \right) \},
\]
 where
$\mathcal{C}_{ol}\left(G \left( \bvsig\right) \right) $ is a
space spanned by the columns of $G\left( \bvsig\right) $.
 Therefore, substituting this solution
into the total complementary function $\Xi$, the canonical dual
problem can be formulated as
 \eb
 (\calP^d): \; \;\ext \left\{\PP^d(\bvsig) = -\half \la  [G(\bvsig)]^{-1}F(\bvsig),  F(\bvsig)\ra
 -V^*(\bvsig)
 \;|\;\bvsig\in \calS_a \right\}. \label{dual}
 \ee
The following theorem was originally presented in general nonconvex systems \cite{GaoBook}.

\begin{theorem}[Analytical Solution and Complementary-dual principle]\hfill \newline
Problem ($\calP^d$) is canonically dual to ($\calP$) in the sense
that if $ \barvsig$ is a critical point of $(\calP^d)$, then
\eb
\bx = [G(\bvsig)]^{-1} F(\bvsig)
\ee
is a critical point of $(\calP)$, the pair $(\barbx, \barvsig)$ is a critical point of $\Xi(\bx, \bvsig)$,   and
 \eb \PP(\barbx) =
\Xi(\barbx, \barvsig)  = \PPd( \barvsig). \label{Equal}\ee
\end{theorem}

This theorem shows that there is no duality gap between the primal
problem ($\calP$) and its canonical dual  ($\calP^d$).
Actually, in $\Xi(\bx, \bvsig)$ the first term
\eb
\Gap(\bx, \bvsig) =  \half \la  \bx ,  G(\bvsig) \bx \ra
\ee
is the  { complementary gap function}, first
introduced  by Gao and Strang in 1989 \cite{GaoStrang89}.
They  proved that if this gap function is positive,  the critical point
$\barbvsig$ is a global maximizer of $\Pi^d$ and the associated $\barbx (\barbvsig) $ is a global minimizer
of the primal problem $(\calP)$.
By introducing the following
notations
 \begin{eqnarray}
 \calS_a^+ &=& \left\{ \bvsig\in
\calS_a\;|\;G\left( \bvsig\right)  \succeq
0\right\}  \label{SaPos},  \\
  \calS_a^- &=& \left\{ \bvsig\in
 \calS_a\;|\;G\left( \bvsig\right) \prec 0\right\}
\label{Sa-} ,
 \end{eqnarray}
where $\mathbf{G}\left( \bvsig\right)  \succeq 0$ means that
$\mathbf{G}\left( \bvsig\right)$ is positive semi-definite and
$\mathbf{G}\left( \bvsig\right) \prec 0$ means that
$\mathbf{G}\left( \bvsig\right)$ is negative definite,
the  Gao and Strang  canonical min-max duality theory can be stated as
\eb
\Pi(\barbx) = \min_{\bx \in \real^n}  \Pi(\bx)  = \max_{\bvsig \in \calS^+_a }  \Pi^d(\bvsig) = \Pi^d(\barbvsig) . \label{eq-minmax}
\ee
This general result has been used extensively in nonconvex analysis and mechanics \cite{GaoBook,yau-gao}.
In 1996,  it was discovered by Gao
that if the gap function is negative in a neighborhood $\calX_o \times \calS_o \subset
\real^n \times \calS^-_a $
of $(\barbx, \barbvsig)$, then    either
the double-max duality relation
\eb
\Pi(\barbx) = \max_{\bx \in \calX_o}  \Pi(\bx)  = \max_{\bvsig \in \calS_o}
 \Pi^d(\bvsig) = \Pi^d(\barbvsig) \label{eq-bimax}
\ee
holds or the double-min duality relation
\eb
\Pi(\barbx) = \min_{\bx \in \calX_o} \Pi(\bx)  = \min_{\bvsig \in \calS_o}   \Pi^d(\bvsig) = \Pi^d(\barbvsig). \label{eq-bimin}
\ee
 Therefore, the   { triality theorem}  was formed by these three pairs of dualities and
  has been used extensively
 in nonconvex mechanics \cite{GaoBook,gao-ogden} and
 global optimization \cite{fang-gaoetal07,Gao-Sherali-AMMA09,r-g-j}.
However, it was realized  in 2003
\cite{gao-opt03,Gao-amma03} that  if the dimensions of the primal   problem
and its canonical dual
 are different, the double-min duality (\ref{eq-bimin}) needs  ``certain
additional conditions".
For the sake of mathematical rigor, the double-min duality was not included in the triality
theory and these additional constraints were left as an open problem (see Remark 1 in \cite{gao-opt03}, also
Theorem 3 and its Remark in a review article by Gao \cite{Gao-amma03}).
By the facts that the double-max duality (\ref{eq-bimax}) is always true and  the double-min duality
plays a key role in  real-life  applications, it was still included in the triality theory
in the "either-or" form  in many applications  for the purposes of
perfection in esthesis and some other reasons in reality.
In the following sections, we will show that the triality theorem holds strongly for the problems
 it was originally proposed. Also we will explain   the reasons
 why  the ``certain additional conditions" in the double-min duality were ignored.

\section{Strong Triality Theory}
In the case  $n=m,$ the triality theorem holds strongly in the following  form.

\begin{theorem}[Tri-duality Theorem]
Suppose that $\barvsig$
 is a critical point of the canonical problem
$ ( \calP^d) $ and  $\barbx =\left[ G \left( \barvsig\right)
\right] ^{-1} F(\barvsig)$.

If $\barvsig\in \mathcal{S}_{a}^{+},$ then $\barvsig$ is a global
maximizer of Problem $ ( \mathcal{P} ^{d} )$ in
$\mathcal{S}_{a}^{+}$ if and only if $\barbx$ is a global
minimizer of Problem $\left( \mathcal{P} \right)$, i.e.,
 the following canonical min-max duality statement holds:
\begin{equation}
\PP (\barbx)=\min_{\bx\in \mathbb{R} ^{n}}\PP \left( \bx\right)
\Longleftrightarrow  \max_{ \bvsig\in \mathcal{S}_{a}^{+}}\PP
^{d}\left( \bvsig \right) =  \PP ^{d}(\barvsig) . \label{Global}
\end{equation}

If $\barvsig\in \mathcal{S}_{a}^{-}$, then  there exists  a neighborhood
$\mathcal{X}_{o}\times \mathcal{S}_{o}\subset \mathbb{R}
^{n}\times \mathcal{S}_{a}^{-}$ of $\left( \barbx,\barvsig\right)
$ such that   we have either  the double-min duality statement
\begin{equation}
\PP (\barbx)=\min_{\bx\in \mathcal{X}_{o}}\PP \left( \bx\right) \;
\Longleftrightarrow  \;  \min_{ \bvsig\in \mathcal{S}_{o}}\PP
^{d}\left( \bvsig\right) = \PP ^{d}\left( \barvsig\right) ,
\label{double-min}
\end{equation}
 or the double-max duality statement
\begin{equation}
\PP (\barbx) = \max_{\bx\in \mathcal{X}_{o}}\PP \left( \bx\right)
\; \Longleftrightarrow  \; \max_{ \bvsig\in \mathcal{S}_{o}}\PP
^{d}\left( \bvsig\right)  = \PP ^{d}\left( \barvsig\right) .
\label{double-max}
\end{equation}
\end{theorem}

\noindent \textbf{Proof.}
If  $(\barbx, \barvsig)$
 is a critical point of the total complementary function $\Xi(\bx, \bvsig)$,
 then by Theorem 1, we have   $\barbx =\left[ G \left( \barvsig\right)
\right] ^{-1} F(\barvsig)$, and
 \eb
 \nabla^2 \PP^d(\barvsig) = -(\nabla\Lbd(\barbx))^T
 [G(\barvsig)]^{-1} \nabla\Lbd(\barbx) -  \nabla^2
 V^*(\barvsig),  \label{HesDual}
 \ee
 \eb
  \nabla^2 \PP(\barbx) = G(\barvsig) + \nabla\Lbd(\barbx)
  \nabla^2V(\Lbd(\barbx)) (\nabla\Lbd(\barbx))^T. \label{HesP}
 \ee
By  the assumption
(A2) we know that $V(\bxi)$ is strictly convex, then,
 \eb
 \nabla^2(V(\Lbd(\barbx))) = (\nabla^2 V^*(\barvsig))^{-1} \succ
 0, \label{Twe}
 \ee
 where $\bar{\bxi} = \Lbd(\barbx)$. Substituting (\ref{Twe}) into (\ref{HesP}),
 we obtain
 \eb
  \nabla^2 \PP(\barbx) = G(\barvsig) + \nabla\Lbd(\barbx)
   (\nabla^2 V^*(\barvsig))^{-1} (\nabla\Lbd(\barbx))^T.
  \label{HesPrimal}
 \ee
\begin{itemize}
    \item Proof of the canonical min-max duality statement (\ref{Global})
    (this proof is a finite-dimensional version
of Gao and Strang's original proof of Theorem 2 in nonconvex analysis \cite{GaoStrang89}).
\end{itemize}
Suppose that $\barbvsig \in \calS^+_a$  is a critical point.
  Since $\Pi^d(\bvsig)$  is concave on $\calS^+_a$, the critical point
 $\barbvsig \in \calS^+_a$  must be a global maximizer of $\Pi^d(\bvsig) $ on $\calS^+_a$.

 On the other hand, if $\barvsig\in \mathcal{S}_{a}^{+}$, the gap function
$\Gap(\bx, \barbvsig) = \half \la  \bx ,  G(\barvsig)
\bx \ra $ is  convex in $\bx \in \real^n$. By the convexity of  $V:\calV_a \rightarrow \real$, we have
\cite{GaoStrang89}
 \begin{eqnarray*}
  \PP(\bx) - \PP(\barbx)
 &  \geq  &  \la \nabla V(\Lbd(\barbx)) ;  \Lbd(\bx)
 -\Lbd(\barbx)\ra +  \half \la \bx ,
 A\bx \ra - \half  \la \barbx ,  A\barbx \ra - \la \bx - \barbx, \bff  \ra \\
 &=&  \Gap(\bx, \barvsig)   -  \Gap(\barbx, \barvsig)  - \la  \bx - \barbx, F(\barvsig)\ra \\
 &\geq &  \la \bx - \barbx, G(\barvsig)\barbx - F(\barvsig) \ra =0 \;\; \forall \bx \in \real^n.
 \end{eqnarray*}
Thus, $\barbx = [G(\barvsig)]^{-1}F(\barvsig)$ is a global
minimizer of problem ($\calP$). Furthermore, $\barvsig$ is also a
global maximizer of Problem $ ( \mathcal{P} ^{d} )$ in
$\mathcal{S}_{a}^{+}$ and the statement (\ref{Global}) holds by Theorem 1.
\begin{itemize}
    \item Proof of the double-min duality statement
    (\ref{double-min}).
\end{itemize}
Suppose that $\barvsig\in\calS_a^-$ and $\barvsig$ is a local
minimizer of problem ($\calP^d$). Then, we have $
 \nabla^2\PP^d(\barvsig)  \succeq 0 $
 and
 \eb
 -(\nabla\Lbd(\barbx))^T
 [G(\barvsig)]^{-1} \nabla\Lbd(\barbx) \succeq  \nabla^2
 V^*(\barvsig) \succ 0. \notag
 \ee
Thus, $\nabla\Lbd(\barbx)$ is invertible, which leads to
 \eb
- G(\barvsig) \preceq \nabla\Lbd(\barbx) ( \nabla^2
 V^*(\barvsig))^{-1} (\nabla\Lbd(\barbx))^T.
 \ee
Therefore, we have
 \eb
  \nabla^2 \PP(\barbx) = G(\barvsig) + \nabla\Lbd(\barbx)
   (\nabla^2 V^*(\barvsig))^{-1} (\nabla\Lbd(\barbx))^T \succ  0.
 \ee
By the assumption (A3), $\barbx = [G(\barvsig)]^{-1}F(\barvsig)$
is also a local minimizer of problem ($\calP$). The reversed  statement can be proved
in the similar way. Thus, (\ref{double-min}) holds.
\begin{itemize}
    \item Proof of the double-max duality statement
    (\ref{double-max}).
\end{itemize}
Suppose that $\barvsig\in \mathcal{S}_{a}^{-}$ and $\barvsig$ is a
local maximizer of problem ($\calP^d$). Then,
$\nabla^2\PP^d(\barvsig)\preceq 0.$ By Theorem 1, $\barbx
= [G(\barvsig)]^{-1}F(\barvsig)$ is a critical point of problem
($\calP$). Due to the assumption (A3), $\nabla^2\PP(\barbx)$ is
invertible. By the well-known Sherman-Morrison-Woodbury identity \cite{Matrix},
 $\nabla^2\PP^d(\barvsig)$ is also
invertible. Furthermore,
\begin{eqnarray*}
(\nabla^2\PP(\barbx))^{-1} = G(\barvsig)^{-1} + G(\barvsig)^{-1}
\nabla \Lbd(\barbx) (\nabla^2\PP^d(\barvsig))^{-1} (\nabla
\Lbd(\barbx))^{T} G(\barvsig)^{-1}\prec 0.
\end{eqnarray*}
Thus, $\barbx = [G(\barvsig)]^{-1}F(\barvsig)$ is also a local
maximizer of problem ($\calP$). Similarly, we can prove the reversed statement.
 Therefore, the triality theorem holds strongly for the
case $n=m$. \hfill $\Box$

\begin{remark}
The tri-duality
  theorem provides
global extremum criteria for three types solutions of the
nonconvex problem $(\calP)$:
  a  global minimizer $\barbx(\barbvsig)$ if   $\barbvsig \in \calS^+_a$
  and a pair of the largest-valued local extrema, i.e.,
   $\barbx(\barbvsig)$  is a global maximizer (resp. minimizer)  if $\barbvsig \in \calS^-_a$ is
a  local maximizer (reps. minimizer).
This pair of largest local extrema plays a critical  role in nonconvex
mechanics  and phase transitions.
\end{remark}
\begin{remark}
The tri-duality theorem  can also be used to identify saddle points of the
primal problem, i.e.  $\barbvsig \in \calS^-_a$ is
a   saddle point of $\Pi(\bvsig)$  if and only if  $\barbx =  [G(\barbvsig)]^{-1}F(\barbvsig)$
 is a
saddle point of $\Pi(\bx)$.
By the facts that  the saddle points are not stable and
 do not exist physically, these points are excluded
from the triality theory.
\end{remark}

The triality theory was first discovered in
post-buckling analysis of a large deformed elastic beam model proposed by Gao
in 1996 \cite{gao-tri96,gao-amr}, where the primal functional is a double-well potential
of  a 2-dimensional displacement field, and its canonical dual is the
so-called {\em pure complementary energy} defined on a 2-dimensional stress field.
Therefore, the triality theory was first proposed  in its strong form,
 i.e. the  tri-duality theorem.

\section{Triality Theory for General Case}

We now consider  the general case   $m\neq n$.
Suppose that $\barbx$ and $\barvsig$ are the
critical points of problem ($\calP$) and ($\calP^d$),
respectively, where $\barbx=[G(\barvsig)]^{-1}F(\barvsig)$ and
$G(\barvsig)$ is invertible. In this case, we also can show that
\eb
  \nabla^2 \PP(\barbx) = G(\barvsig) + \nabla\Lbd(\barbx)
   (\nabla^2 V^*(\barvsig))^{-1} (\nabla\Lbd(\barbx))^T,
  \label{HesPri}
 \ee
 and
\eb
 \nabla^2 \PP^d(\barvsig) = -(\nabla\Lbd(\barbx))^T
 [G(\barvsig)]^{-1} \nabla\Lbd(\barbx) -  \nabla^2
 V^*(\barvsig). \label{HesDua}
 \ee

Suppose that $m<n$. By the  Sherman-Morison-Woodbury Theorem
in \cite{Matrix} and the assumption (A3), we have
 \eb
 [\nabla^2 \PP(\barbx)]^{-1} =
 [G(\barvsig)]^{-1}+[G(\barvsig)]^{-1}\nabla\Lbd(\barbx) (\nabla^2
 \PP^d(\barvsig))^{-1} (\nabla\Lbd(\barbx))^T
   [G(\barvsig)]^{-1}.
 \ee
This shows that $\nabla^2 \PP^d(\barvsig)$ is  invertible. Similarly,
 we can show that $\nabla^2 \PP^d(\barvsig)$ is  also
invertible if  $m>n$.

\begin{lemma}
Suppose that $m<n$, the critical point $\barvsig\in\calS_a^-$ is a
local minimizer of problem ($\calP^d$). Then, $\nabla^2
\PP(\barbx)$ has $m$ positive eigenvalues and $n-m$ negative
eigenvalues, i.e., there exists two matrices
$\Pf\in\R^{n \times m}$ and $ \Ps\in\R^{n\times (n-m)}$ such that
 \eb
 \Pf^T \nabla^2 \PP(\barbx) \Pf \succ 0\;{\rm and}\;
 \Ps^T \nabla^2 \PP(\barbx) \Ps \prec 0.
 \ee
\end{lemma}
\noindent \textbf{Proof.} Since $\nabla^2 V^*(\barvsig)\succ 0$,
there exists a invertible matrix $R\in\R^{m\times m}$ such that
$\nabla^2 V^*(\barvsig)=R^T R$. Thus, we have
 \eb
 - (\nabla\Lbd(\barbx) R^{-1})^T[G(\barvsig)]^{-1} \nabla\Lbd(\barbx)
 R^{-1} - I_{m\times m} \succ 0.
 \ee
 Note that $G(\barvsig)\prec 0$ and $\nabla\Lbd(\barbx)
 R^{-1} (\nabla\Lbd(\barbx)
 R^{-1}))^T\succeq 0$. There exists a matrix $T$ such that
  \eb
  T^T G(\barvsig) T = \Diag(-\lambda_1,\cdots,-\lambda_n),
  \ee
  and
  \eb
   T^T \nabla\Lbd(\barbx) R^{-1} (\nabla\Lbd(\barbx) R^{-1}))^T T
 =\Diag (a_1,\cdots,a_m,0,\cdots,0), \label{35}
  \ee
  where $\lambda_k >0,\;k=1,\cdots,n,$ and $a_k>0, k=1,\cdots,m$.
According to the  decomposition theory of singular matrices, we know that
there exist orthogonal matrices $U\in\R^{n\times n}$ and $E\in
\R^{m\times m}$ such that
 \eb
 T^T \nabla\Lbd(\barbx) R^{-1} = U \left(
\begin{array}{ccc}
  \sqrt{a_1} &  &  \\
   & \ddots &  \\
  &  & \sqrt{a_m} \\
   0 & \cdots & 0  \\
   & \cdots &  \\
    0 & \cdots & 0  \\
\end{array}
\right) E.
 \ee
In light of (\ref{35}), we know that $U =I_{n\times n}$. Then, we
have
\begin{eqnarray*}
 (R^{-1})^T \nabla^2\PP^d(\barvsig) R^{-1} &=&  - (\nabla\Lbd(\barbx) R^{-1})^T[G(\barvsig)]^{-1} \nabla\Lbd(\barbx)
 R^{-1} - I_{m\times m}\\
 &=& - (T^T \nabla\Lbd(\barbx) R^{-1})^T[T ^T G(\barvsig)T]^{-1}
T \nabla\Lbd(\barbx)
 R^{-1} - I_{m\times m}\\
 &=& E^T \Diag \left(\frac{a_1}{\lambda_1}-1,\cdots,\frac{a_m}{\lambda_m}-1
 \right)E
 \succ 0.
\end{eqnarray*}
Thus, $a_k>\lambda_k,\;k=1,\cdots,m$. It is easy to verify that
 \eb
 T^T \nabla^2\PP(\barbx)T = \Diag(a_1-\lambda_1,\cdots,
 a_m-\lambda_m,-\lambda_{m+1},\cdots,-\lambda_n).
 \ee
This shows that  $\nabla^2 \PP(\barbx)$ has $m$ positive eigenvalues and
$n-m$ negative eigenvalues. Therefore, the matrix  $\Pf$ can be obtained by
collecting all the eigenvectors corresponding to the positive
eigenvalues and $\Ps$ can be obtained by collecting all the
eigenvectors corresponding to the negative eigenvalues.   \hfill $\Box$

In a similar way, we can prove the following lemma.

\begin{lemma}
Suppose that $m>n$ and the critical point $\barbx=\left[ G\left( \barvsig\right)
\right] ^{-1}F(\barvsig)$ is
a local minimizer of Problem $(\calP)$, where $\barvsig \in
\mathcal{S}_{a}^{-}$. Then, $\nabla^2 \PP^d(\barvsig)$ has $n$
positive eigenvalues and $m-n$ negative eigenvalues, i.e.,
there exists two matrices $\Qf \in\R^{m \times n}$ and $
\Qs \in\R^{m\times (m-n)}$ such that
 \eb
\Qf^T \nabla^2 \PP^d(\barvsig) \Qf \succ 0\;{\rm and}\;
 \Qs^T \nabla^2 \PP^d(\barvsig) \Qs \prec 0.
 \ee
\end{lemma}

Let the $m$ column vectors of $\Pf$ be $\bp^\flat_{1}, \cdots, \bp^\flat_{m}$
and  the $n$ column vectors of $\Qf$ be $\bq^\flat_{1}, \cdots,
\bq^\flat_{n}$, respectively. Clearly, $\bp^\flat_{1}, \cdots, \bp^\flat_{m}$ and
$\bq^\flat_1, \cdots, \bq^\flat_n$ are two sets of linearly independent
vectors, respectively. By introducing two subspaces
 \begin{eqnarray}
  \mathcal{X}_{\flat} &=& \{\bx \in \Rn\;|\;\bx = \barbx
  + \theta_{1}\bp^\flat_{1}+\cdots+\theta_{m}\bp^\flat_{m},\;
  \theta_{i} \in \R,\; i=1,\cdots,m
  \}, \\
  \mathcal{S}_{\flat} & = & \{\bvsig \in \mathbb{R}^{m}\;|\;\bvsig =
  \barvsig + \vartheta_{1}\bq^\flat_{1}+\cdots+\vartheta_{n}\bq^\flat_{n},\;
  \vartheta_{i} \in \mathbb{R},\; i=1,\cdots,n
  \},
  \end{eqnarray}
 the  triality theory holds  for general case in the following refined form.

\begin{theorem}[Triality Theorem] \label{WTT}~~~

Suppose that $\barvsig$
 is a critical point of
problem ($\calP^d$) and $\barbx = [G(\barvsig)] ^{-1}F(\barvsig)$.

If $\barvsig\in \calS_a^+,$ then the canonical min-max duality
holds in the strong form of
\begin{equation}
\PP(\barbx) = \min_{\bx\in \mathbb{R} ^{n}}\PP\left(
\bx\right)\Leftrightarrow  \max_{ \bvsig\in
\calS_{a}^{+}}\PP^{d}\left( \bvsig \right) = \PP^{d}(\barvsig) .
\label{PmGlobal}
\end{equation}

If $\barvsig\in \calS_{a}^{-}$, then there exists a neighborhood
$\mathcal{X}_{o}\times \mathcal{S}_{o}\subset \Rn\times
\mathcal{S}_{a}^{-}$ of $\left( \barbx,\barvsig\right) $ such that
the double-max duality holds in the strong form of
\begin{equation}
\PP(\barbx) = \max_{\bx\in \mathcal{X}_{o}}\PP\left( \bx\right)
\Leftrightarrow \max_{ \bvsig\in \mathcal{S}_{o}}\PP^{d}\left(
\bvsig\right) = \PP^{d}\left( \barvsig\right) .  \label{Pmdmax}
\end{equation}
However,  the double-min duality statement holds conditionally in
the following   symmetrical forms.
\begin{enumerate}
\item   If $m<n$ and $\barvsig\in \mathcal{S}_{a}^{-}$ is a local
minimizer of  $\PPd(\bvsig)$,  then $\barbx= [G(\barvsig)]^{-1}
F(\barvsig)$ is a saddle point of $\PP(\bx)$  and the double-min
duality holds weakly on $ \mathcal{X}_{o}\cap\mathcal{X}_{\flat}
\times \calS_o$, i.e.
 \eb
  \PP(\barbx) = \min_{\bx\in \mathcal{X}_{o}\cap\mathcal{X}_{\flat}}\PP\left( \bx\right)
=\min_{ \bvsig\in \mathcal{S}_{o}}\PP^{d}\left(
\bvsig\right)=\PP^{d}(\barvsig). \label{WDM1}
 \ee

\item If $m>n$ and $\barbx=[G(\barvsig)] ^{-1}F(\barvsig)$ is a
local minimizer of $\PP(\bx),$
 then $\barvsig$ is a saddle point of $\PPd(\bvsig) $
 and  the double-min duality holds weakly on $\mathcal{X}_o \times  \mathcal{S}_{o}\cap
\mathcal{S}_{\flat}$, i.e.
 \eb
  \PP(\barbx) = \min_{\bx\in \mathcal{X}_{o}}\PP\left( \bx\right)
=\min_{ \bvsig\in \calS_{o}\cap \mathcal{S}_{\flat}}\PP^{d}\left(
\bvsig\right)=\PP^{d}(\barvsig).\label{WDM2}
 \ee

\end{enumerate}
\end{theorem}

\noindent \textbf{Proof.} The proof of min-max duality statement
(\ref{PmGlobal}) and the double-max duality statement
(\ref{Pmdmax}) are the same to the proof of (\ref{Global}) and
(\ref{double-max}). We only need to prove (\ref{WDM1}) and
(\ref{WDM2}).

Suppose that $m<n$ and $\barvsig\in \calS_a^-$ is not only a local
minimizer, but also a critical point of problem ($\calP^d$). By Lemma 1
 we know that $\nabla^2\PP(\barbx)$ has both positive and
negative eigenvalues. Thus, $\barbx=[G(\barvsig)]^{-1}F(\barvsig)$
is a saddle point of problem ($\calP$). We  let
 \eb
 \varphi(t_{1},\cdots,t_{m}) =
 \PP(\barbx+t_{1}\bp^\flat_{1}+\cdots+t_{m}\bp^\flat_{m}), \label{phit}
 \ee
where $\bp^\flat_1,\cdots,\bp^\flat_m$ are the column vectors of $\Pf$ defined in
Lemma 1. By direct verification, we have
 \eb
 \nabla \varphi(0,\cdots,0) = (\nabla\PP(\barbx))^T \Pf =0
 \ee
 and
  \eb
  \nabla^2 \varphi(0,\cdots,0) = \Pf^T \nabla^2\PP(\barbx) \Pf
  \succ 0.
  \ee
 Thus, $(0,\cdots,0)$ is a local minimizer of
 $\varphi(t_1,\cdots,t_m)$. Hence, the equation (\ref{WDM1}) holds.  The statement
 (\ref{WDM2}) can be proved in the similar way.  \hfill
 $\Box$\\

\begin{remark}[NP-hard Problems and Perturbation]
The canonical min-max duality (\ref{PmGlobal}) shows that the nonconvex minimization problem is
 equivalent to a concave maximization dual problem over a closed convex set $\calS^+_a$.
 If $\Pi^d(\bvsig)$ has at least one critical point in $\calS^+_a$,
 the global minimizer of $\Pi(\bx)$ can be easily obtained by the canonical duality theory.
 However, if $\Pi^d(\bvsig)$ has no critical points in $\calS^+_a$,
 to find global minimizer for nonconvex function $\Pi(\bx)$ could be very difficult.
 If the  vector  $\bff ={\bf 0} \in \real^n$, the problem $(\calP)$ is { homogenous}.
 Moreover, if  $\bb_k = {\bf 0 } \in \real^n  \; \forall k = 1, \cdots, m$,
then the geometrical operator $\Lam(\bx)$ is a pure quadratic measure
 (i.e. an objective  measure in certain space).
 In this case, the vector $F(\bvsig) = {\bf 0}$, the set $\calS^+_a$ is empty, and the canonical dual function
 $\Pi^d(\bvsig)  =  - \VV^*(\bvsig)$ is   concave,  which has only a unique  maximizer $\barbvsig$. By the double-max duality we know that the corresponding primal solution
 $\barbx = {\bf 0} $ is a local maximizer if $\barbvsig \in \calS^-_a$.
 From the point view of systems theory, the pure quadratic operator $\Lam(\bx)$ means  that  the system possesses   certain symmetry.
 If there is no input ($\bff = {\bf 0}$), the primal function $\Pi(\bx)$
 could have multiple global minimizers.
 It was indicated in \cite{gao-jimo07} that a nonconvex minimization problem
 could be NP-hard if its canonical dual has no KKT (or critical) point in $
\calS^+_a$.
In order to solve this type problems,
 several perturbation methods have been suggested in
 \cite{gao-ruan-jogo10,r-g-j,wang-etal}.
 It is shown very recently that by the canonical duality theory, a class of  NP-hard box/integer
 constrained programming problems
 are equivalent to unconstrained canonical dual problems in continuous space, which can be solved via
 deterministic methods \cite{gao-watsonetal}.
 \end{remark}

Dual to $\calX_\flat$ and $\calS_\flat$,  we can let
 \begin{eqnarray*}
\calXs  &=& \{\bx \in \Rn\;|\;\bx = \barbx
  + \theta_{1}\bp^\sharp_{1}+\cdots+\theta_{n-m}\bp^\sharp_{n-m},\;
  \theta_{i} \in \R,\; i=1,\cdots,n-m
  \}, \\
\calSs & = & \{\bvsig \in \mathbb{R}^{m}\;|\;\bvsig =
  \barvsig + \vartheta_{1}\bq^\sharp_{1}+\cdots+\vartheta_{m-n}\bq^\sharp_{m-n},\;
  \vartheta_{i} \in \mathbb{R},\; i=1,\cdots,m-n
  \},
  \end{eqnarray*}
where $\{\bp^\sharp_i \} $ and $\{\bq^\sharp_i \}$ are column vectors  of $\Ps$ and $\Qs$,
 respectively. Then, complementary to the weak double-min duality statements (\ref{WDM1}) and
(\ref{WDM2}), we have the following weak saddle min-max duality theorem.
\begin{theorem}[Weak Saddle Duality Theorem] \label{thm-wsmm}~~~

Suppose that $\barvsig\in \calS_{a}^{-}$
 is a critical point of
problem ($\calP^d$), the vector  $\barbx = [G(\barvsig)] ^{-1}F(\barvsig)$,
and  $\mathcal{X}_{o}\times \mathcal{S}_{o}\subset \Rn\times
\mathcal{S}_{a}^{-}$  is a  neighborhood
of $(\barbx, \barbvsig)$.
\begin{enumerate}
\item   If $m<n$ and $\barvsig\in \mathcal{S}_{a}^{-}$ is a local
minimizer of  $\PPd(\bvsig)$,  then $\barbx= [G(\barvsig)]^{-1}
F(\barvsig)$ is a saddle point of $\PP(\bx)$  and the saddle max-min
duality holds weakly on $ \mathcal{X}_{o}\cap\calXs
\times \calS_o$, i.e.
 \eb
  \PP(\barbx) = \max_{\bx\in \mathcal{X}_{o}\cap\calXs}\PP\left( \bx\right)
=\min_{ \bvsig\in \mathcal{S}_{o}}\PP^{d}\left(
\bvsig\right)=\PP^{d}(\barvsig).
 \ee

\item If $m>n$ and $\barbx=[G(\barvsig)] ^{-1}F(\barvsig)$ is a
local minimizer of $\PP(\bx),$
 then $\barvsig$ is a saddle point of $\PPd(\bvsig) $
 and  the saddle min-max duality holds weakly on $\mathcal{X}_o \times  \mathcal{S}_{o}\cap
\calSs$, i.e.
 \eb
  \PP(\barbx) = \min_{\bx\in \mathcal{X}_{o}}\PP\left( \bx\right)
=\max_{ \bvsig\in \calS_{o}\cap \calSs}\PP^{d}\left(
\bvsig\right)=\PP^{d}(\barvsig).
 \ee

\end{enumerate}
\end{theorem}

  \begin{remark}
Theorem 3 shows that both the canonical min-max and double-max duality statements hold strongly for general cases; the double-min duality holds strongly for n = m but
weakly for $n\neq  m$  in a symmetrical form. The  ``certain additional conditions"
are simply the intersection $\calX_o \bigcap \calX_\flat$ for $n > m$ and
$\calS_o \bigcap \calS_\flat$ for $ n < m$.
 Therefore, the open problem on the double-min duality left in
2003 \cite{gao-opt03,Gao-amma03} is now solved!
While from Theorem 4 we know  that if $G(\barbvsig) \prec 0$ and $n \ge m$,
the solution  $\barbx (\barbvsig)$
 could be a saddle point.
 Mathematically speaking,  nonstable critical points   do
 not produce any computational difficulties in
numerical optimization. Also,  in real-life problems the saddle point  is  not
 considered as a phase state and   does not physically exist.
 These are the part of  reasons why the saddle point in $\calS^-_a$ is  ignored by the triality theory.
\end{remark}

The triality theory has been challenged recently by a large number of counter-examples
in a series of more than seven  papers (see \cite{vz-dcds,vz1} and references cited therein).
It was written  in \cite{vz-dcds} that  ``Because our counter-examples are very simple, using quadratic functions defined
on whole Hilbert (even finite dimensional) spaces, it is difficult to reinforce the
hypotheses of the above mentioned results in order to keep the same conclusions
and not obtain trivialities."
It turns out that in addition to many conceptual mistakes (see Section 6), most of
these counter-examples
 simply  discuss the saddle points in $\calS^-_a$ for the case  $n\neq m$.
In fact, these count-examples
address  the same type of open problem for the double-min duality
left unaddressed  in \cite{gao-opt03,Gao-amma03}\footnote{It is interesting to note  that
the references  \cite{gao-opt03,Gao-amma03} never been  cited in any one of this set of papers.
}.
Indeed, by Theorem 3 we know that the  double-min duality holds conditionally when
$n\neq m$.
Based on Theorems 2 and  3,  we know that   the saddle points could exist in $\calS^-_a$
even if $n = m$; While by Theorem 4
  one can easily construct many other {\em V-Z type counter-examples} which are physically useless.

 \section{Application}

Let us consider the following quadratic-log optimization problem:
 \eb
 (\calP):\; \ext \left\{ \PP(\bx)=\half \bx^T A \bx  -\sum_{k=1}^m \log(\half \bx^T B^k \bx +
 d^k)  -\bx^T \bff  \; |  \;\;  \bx\in \Rn \right\} ,
 \ee
where $A$ is a positive definite matrix, $B^k, k=1,\cdots, m$ are
positive semi-definite matrices and $d^k>0, k=1,\cdots, m$. In
this case, its canonical dual problem can be expressed as:
 \eb
 (\calP^d):\;\; \ext \left\{\PP^d(\bvsig) = -\half \bff^T G(\bvsig)^{-1}\bff
 +\sum_{k=1}^m(d^k\varsigma_k +1 +\log(-\varsigma_k))\;|\;\bvsig\in
 \calS_a\right\},
 \ee
where $\calS_a =
\left\{\bvsig= \{ \varsigma_i \} \in \real^m \;|\;- \frac{1}{d^k} \le \varsigma_k < 0
, k=1,\cdots, m \right\}$.

Let $n=2, \;\; m=1$,  and
 \eb
 A = \left(%
\begin{array}{cc}
  1 &  \\
   & 2 \\
\end{array}%
\right), B^1 = \left(%
\begin{array}{cc}
  5 &  \\
   & 4 \\
\end{array}%
\right), \bff = \left(%
\begin{array}{c}
  0.5 \\
  0.1 \\
\end{array}
\right)\; {\rm and}\; d^1 =1. \notag
 \ee
In this case, we have $\calS_a=\{\varsigma\in\R\;|\;-1\leq\varsigma<0 \}$ and
 \eb
 \PP^d(\varsigma) = -\half \left(\frac{0.5^2}{1+5\varsigma}+\frac{0.1^2}{2+4\varsigma}
 +\varsigma+1+\log(-\varsigma)\right). \notag
 \ee

\begin{figure}[tbph]
\begin{center}
\includegraphics[ width=3.5 in, bb=120 278 474 560]
{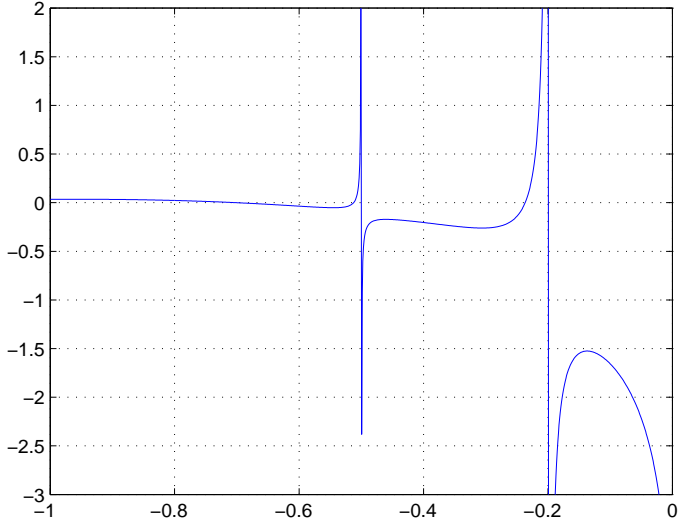}
$\;$ \hspace{4cm} \caption{Grapy of $\PP^d(\varsigma)$}
\end{center}
\end{figure}

\begin{figure*}[tbp]
\centering
\begin{minipage}[b]{0.4\textwidth}
    \centering
    \includegraphics[width=2in]{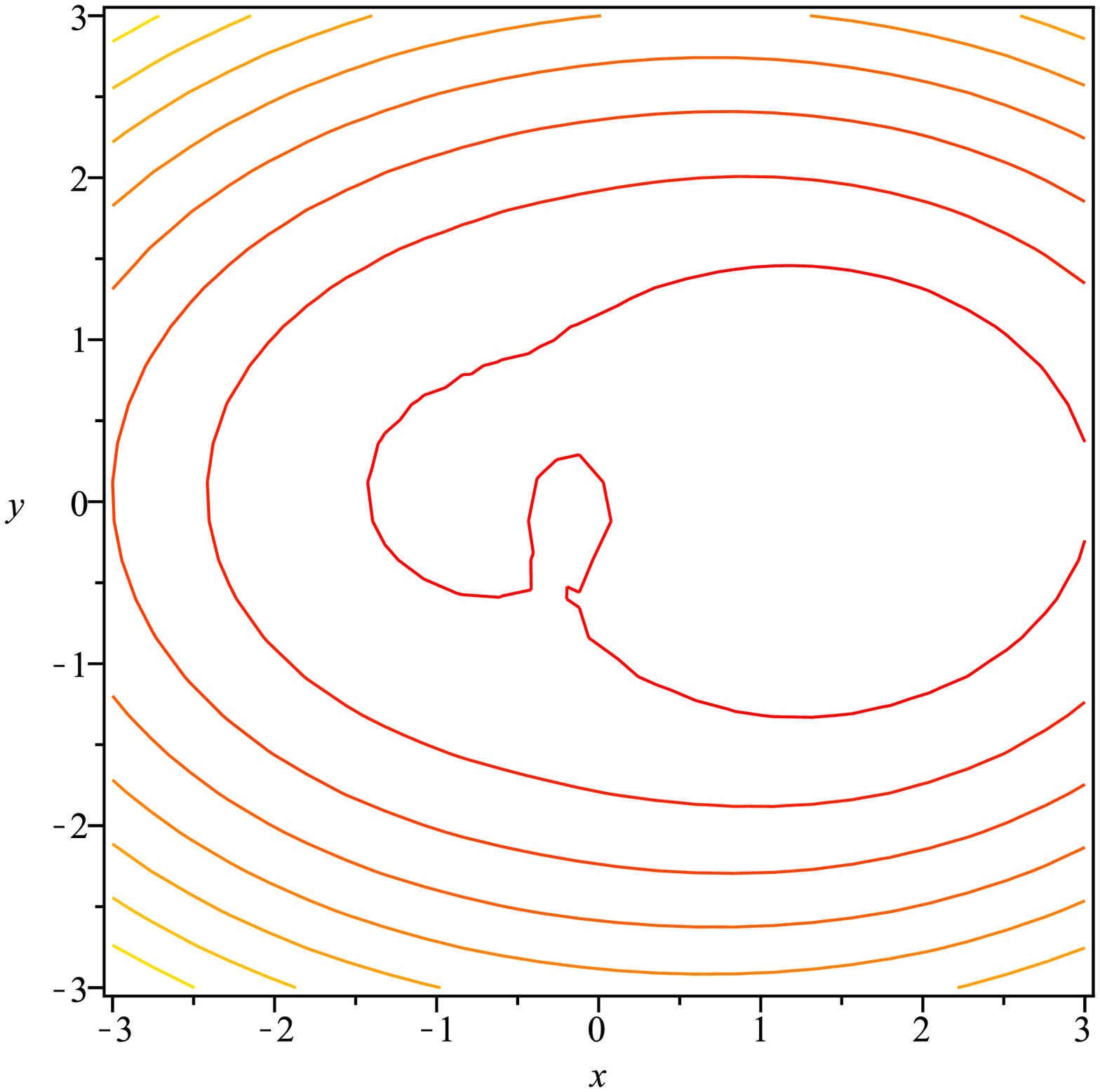}
    \makebox[0.5cm]{\small{a}}
  \end{minipage} \hspace{0.04\textwidth}
\begin{minipage}[b]{0.4\textwidth}
    \centering
    \includegraphics[width=2in]{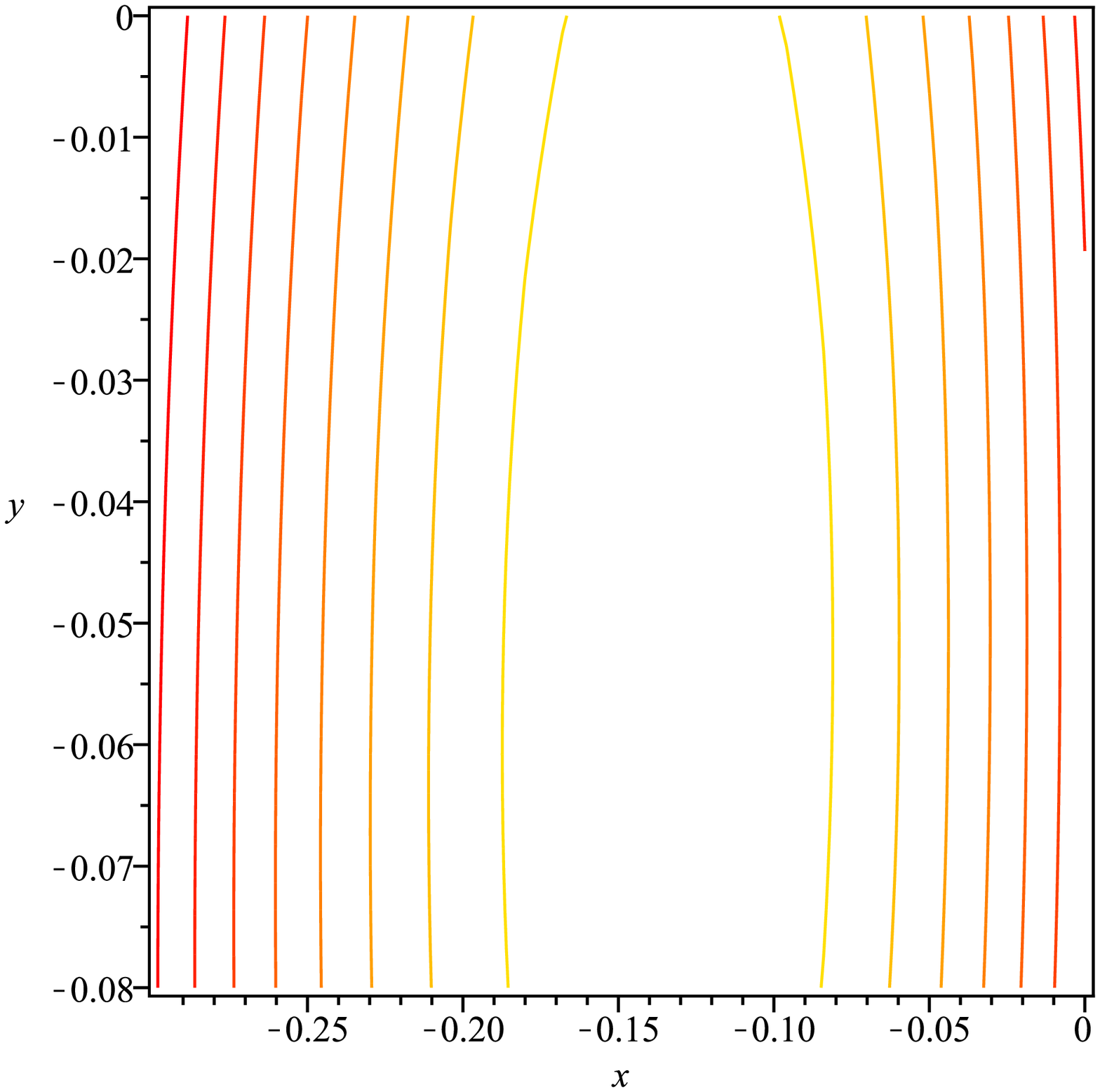}
    \makebox[0.5cm]{\small{b}}
  \end{minipage}\\[5pt]
   \caption{Contours of $\Pi(\bx)$:  (a) global minimizer of $\PP(x)$ and
   (b)  local maximizer  of $\PP(x)$}
\end{figure*}

From  Fig. 1 we can see that
  $\PP^d(\varsigma)$ has one critical point
$\bar{\varsigma}_1 = -0.13696432 $ in $ \calS_a^+ = (-0.2, 0)$
and two critical points
$\bar{\varsigma}_2 = -0.54470504  $ and
$\bar{\varsigma}_3 = -0.95209751 $ in $\calS_a^- = (-1, -0.5)$. Thus, by the triality theorem we know that
$\barbx_1=
(A+\bar{\varsigma}_1 B^1)^{-1}\bff = [1.58640312,
 0.06886375]^T$ is the global minimizer to the problem ($\calP$).
Since  $\bar{\varsigma}_2 \in \calS_a^-$ is a
 local minimizer, $\bar{\varsigma}_3 \in \calS_a^-$ is a
 local maximizer to the problem ($\calP^d$), and  $m=1<n=2$,
 by Theorem 3 we know that $\barbx_2=
(A+\bar{\varsigma}_2 B^1)^{-1}\bff = [-0.2901031,
 -0.5592211]^T$ is a saddle point of $\Pi(\bx)$, while
 $\barbx_3=
(A+\bar{\varsigma}_3 B^1)^{-1}\bff = [ -0.13296148,
 -0.0552978]^T$ is a local maximizer to the problem ($\calP$).
 More applications can be found in \cite{gao-wu-jimo11}.

\section{Some Fundamental Concepts in  Canonical Systems}
 Global optimization problem in mathematics is usually formulated in the following general form
\[
\min \{ f(\bx) | \; \bx \in \calX \subset \real^n \},
\]
where the real-valued function $f(\bx)$ is simply assumed to be nonconvex (or Lipschitz,  differentiable, etc.) on its feasible space $\calX \in \real^n$, in which certain constraints are given.
It is known that  this problem  could have a large number of local extrema
and to identify global optima is a main challenging task in global optimization.
If there is no detailed information  available for the given function $f(\bx)$,
it is  difficult (may be impossible)
to have a general  theory and method for solving
this general problem effectively. Also, to find the largest local extrema is fundamentally important
in many real-life applications\footnote{
It should be  emphasized here that to find the largest local maximum of $f(\bx)$
is not simply equivalent to solve the problem $\min \{ -  f(\bx) | \;\bx \in \calX \}$.}.

Mathematics and physics (mechanics) have been complementary partners since Newton's time.
 It is known that the calculus of variation and mathematical optimization
 were originally developed from Euler-Lagrange mechanics.
Also, the modern mathematical theory of convex analysis was started from J.J. Moreau's pioneering work  in contact mechanics \cite{moreau}.
However,  as V.I. Arnold pointed out \cite{arnold}: ``In the middle of the twentieth century it was attempted to divide physics and mathematics. The consequences turned out to be catastrophic."
For example, in mathematical physics, the objectivity is directly related to some fundamental concepts
and principles, such as
geometrical nonlinearity, constitutive laws,  and  work-conjugate principle, etc.
A function(al) can be  called objective or free energy only if  certain
 intrinsic constraints (physical laws) are satisfied (see \cite{gao-ogden}).
Unfortunately,
  the  objective function in mathematical optimization
has  been misused with other concepts such as cost function, energy function, and energy functional\footnote{See the web page at http://en.wikipedia.org/wiki/Mathematical\_optimization},
which  leads to  some conceptual mistakes.
This section will discuss   some   important issues
in classical Lagrangian mechanics/duality, mathematical optimization,  and general systems theory.

\subsection{Canonical systems}

According to E. Tonti \cite{tonti72}, in virtually every physical system there exists at least
three types of variables:

 (1) the configuration  variable $\bx \in \calX $,
  which describes the state or output of the system, such as
 the Lagrangian generalized coordinates (or displacements)  in analytical mechanics \cite{landau-lif},
 decision variable in game theory, etc.

  (2) the source variable $\bx^*  = \bff \in \calX^* $,
  which represents the input of the system, such as the external force in mechanics and  charge density
  in theory of electrical field, etc.

  (3) a pair of internal (or intermediate) variables
$( \bveps, \bveps^*) \in \calE \times \calE^*  $,
which describes certain interior (constitutive) properties of the system,
such as strain and stress in elasticity,
velocity and momentum in dynamics, etc.

 By the facts  that the constitutive laws
  should be objective (coordinates-free) and physical variables appear always in one-to-one pairs
  (i.e. the  Hill work-conjugacy principle in continuum mechanics \cite{GaoBook}),
  it is reasonable to assume that for a given natural  system,
there exists   a
certain objective measure $\bveps = \barLam(\bx):\calX_a \subset \calX \rightarrow \calE_a \subset \calE$ and
a {\em stored energy} $\barW:\calE_a \rightarrow \real$  such that the constitutive duality relation
$\bveps^*= \nabla \barW(\bveps)  :\calE_a \rightarrow \calE^*_a \subset \calE^*$ is canonical (i.e., one-to-one on $\calE_a \times \calE^*_a$).
Such a system is the so-called
{\em canonical system} and is denoted as
$\mathbb{S}_a = \{ \la \calX_a ,\calX_a^* \ra ,   \la \calE_a ;  \calE^*_a \ra ; \barLam, C \}$
(see Chapter 4, \cite{GaoBook}),
where $C = \nabla \barW: \calE_a \rightarrow \calE^*_a$ represents the constitutive mapping,
$\la * , * \ra$ and $\la * ; * \ra$ denote the bilinear forms on $\calX \times \calX^*$ and
$\calE \times \calE^*$, respectively.
The system is called {\em geometrically nonlinear} (resp. linear) if the geometrical operator
$\barLam$ is nonlinear (resp. linear);
the system is called {\em physically (or constitutively) nonlinear } (resp. linear)
if the constitutive operator $C$ is nonlinear (resp. linear);
the systems is called {\em fully nonlinear} (resp. linear) if
it is  both geometrically and physically nonlinear (resp. linear).

 The most simple  geometrically linear  system is   controlled by the  quadratic function
$
 \Pi(\bx)=   \half \la    \bx , A   \bx \ra - \la \bx, \bff \ra  ,
$
where  $A \in \real^{n\times n} $ is a symmetrical matrix\footnote{
The skew symmetric matrix $A_s =\half (A -A^T) $ does not store energy since  $\bx^T A_s \bx  \equiv  0$.}.
If  $A$ is positive (semi) definite, by
    Cholesky decomposition we know that
 there exists a matrix $ D :  \real^n \rightarrow \real^m$   such that $A= D^T D$.
    Therefore, we have
    $ \half \la \bx, A \bx \ra =  \half \la D \bx ; D \bx \ra = \TT(D \bx)$ and
    $\TT(\by)$ is an objective function of   $\by = D \bx  \in \real^m$.
        By the fact that any symmetrical matrix can be written in difference of two positive definite matrices,
    it turns out that any given  quadratic  function
    can be written in the so-called  d.c. (difference of convex functions) form.

\subsection{Geometrically linear systems and Lagrangian duality}

In fact,  the most popular  Lagrangian in its original form is actually  defined  by
\[
\Pi(\bx) = \TT(D\bx) - \UU(\bx),
\]
where  the  objective function   $\TT(\by):\calY_a \subset \real^m \rightarrow \real$
is a kinetic  energy, while
  $\UU:\calX_a \subset \real^n \rightarrow \real$ is a  potential energy of the system\footnote{The Lagrangian form was first introduced by
W. Hamilton in classical mechanics  and denoted by $L = T - U$, which is the standard notation
extensively used from  dynamical systems to quantum field theory (see \cite{landau-lif}).}, which could be linear or convex
  such that $\Pi(\bx)$ is well-defined on
  the so-called {\em kinetically admissible space} $\calX_k = \{\bx \in \calX_a | \; \DD \bx \in \calY_a\}$ \cite{GaoBook}.
For  Newtonian mechanics,
  $\TT(\by)$ is quadratic and  the objectivity of this kinetic energy ensures the validity of Newton's laws under the Galilean transformation;
while for  Einstein's special relativity theory,
the objective function $\TT(\by)$ is strictly convex (see Chapter 2, \cite{GaoBook}),
which is an invariant under the Lorentz transformation.
In either case, the so-called {\em complementary energy } $\TT^*(\by^*)$ can
be uniquely defined on $\calY_a^* \subset \real^m$
by the classical Legendre transformation
$\TT^*(\by^*) =  \sta  \{ \la \by ;  \by^* \ra - \TT(\by) | \;  \by \in \calY_a \}$
such that the original Lagrangian  $\Pi(\bx)$ is equivalent to its mixed form
\eb
L(\bx, \by^*) = \la D \bx ; \by^* \ra - \TT^*(\by^*) - \UU(\bx) : \calX_a \times \calY^*_a \rightarrow \real
\ee
  which is the standard  form in mathematical optimization.
  For conservative systems, the Lagrangian should be a constants, therefore, the
  criticality condition $\nabla L(\bx, \by^*) = 0$
  leads to the well-known {\em Euler-Lagrange equations}:
  \eb
  \DD^* \by^* = \nabla \UU(\bx) , \;\; \DD \bx   = \nabla \TT^*(\by^*), \label{eq-el}
  \ee
  where $\DD^*$ is an adjoint of the linear operator $\DD$ defined by
  $\la \DD \bx ; \by^* \ra = \la \bx , \DD^* \by^* \ra$.
  By the canonical duality  $\by = \nabla \TT^*(\by^*) \;\; \Leftrightarrow \;\;
  \by^* = \nabla  \TT(\by)$, we have
  the equilibrium equation
  \[
  \DD^* \nabla \TT(\DD \bx) = \nabla \UU(\bx).
  \]
  Particularly, if $\UU(\bx) = \la \bx, \bff \ra$
  is linear and $\TT(\by) = \half \la \by ; C \by \ra$ is quadratic,
   where $C$ is a linear operator,
   the equilibrium equation takes  a particular symmetrical form
$\DD^* C \DD \bx = \bff $ which is repeated throughout the field equations of mathematical physics
 \cite{gs86}.

For geometrically linear  static systems,  both the input and the configuration variables are time independent.
In this case, the convex objective function
 $\TT(\by)$ is the so-called internal (or stored) energy and
   $\UU(\bx)$ is   the external potential,
 which should be linear $\UU(\bx)  = \la \bx, \bff\ra $
 such that its derivative $\nabla \UU(\bx) = \bff$ is a given source of the system.
Therefore, the Lagrangian form $\Pi(\bx)$ represents the {\em total potential}
 of the system, which is convex
 on  $\calX_k$   
  and its mixed form $L(\bx, \by^*)$ is a saddle function on $\calX_a \times \calY^*_a$.
Therefore, the  traditional saddle Lagrangian  duality theory
links the convex primal problem $\min \{ \Pi(\bx)| \; \bx \in \calX_k\}$
 to a unique  dual problem
\eb
\max \left\{ \Pi^*(\by^*) = - \TT^*(\by^*) | \; \;  \by^* \in \calY_s^*  \right\},
\ee
where $\calY^*_s = \{ \by^* \in \calY^*_a | \; \DD^* \by^* = \bff \in \calX^*_a \subset \real^n \}$
is the so-called {\em statically admissible space}.
The objectivity 
of this dual problem
 is guaranteed by the objectivity of  $\TT(\by)$.
By introducing a  Lagrange multiplier $\bx$, which must be a solution to the primal problem
(see Lagrange multiplier's law in Section 1.5 \cite{GaoBook}),
 to relax the
 equilibrium constraint $\DD^* \by^* = \bff$ in $\calY^*_s$,
 the Lagrangian  is exactly the mixed form $L(\bx, \by^*)$
 and
  the  one-to-one Lagrangian saddle min-max duality
  \[
  \min_{\bx \in \calX_k} \Pi(\bx) = \min_{\bx \in\calX_a} \max_{\by^*\in \calY^*_a}  L(\bx,\by^*)= \max_{\by^*\in \calY^*_a}  \min_{\bx \in\calX_a} L(\bx,\by^*)
   =  \max_{\by^* \in \calY^*_s} \Pi^*(\by^*)
  \]
 is  called the {\em mono-duality} in    canonical  systems theory \cite{GaoBook}.
 In mathematical economics,   where the objective function $\TT(\DD \bx)$
 is corresponding to the revenue, denoted by $R(\bx)$, and  the potential $\UU(\bx)$
 is the {\em cost function,}  denoted by $ C(\bx)$, then
  $\Pi(\bx) = R(\bx) - C(\bx)$ is the so-called total profit.
   For geometrically linear static problems, the cost function $C(\bx)$ is usually linear, while the revenue $R(\bx)$
    is  a concave objective  function of certain measure (norm) of $\bx$ in order to have maximum total profit $\Pi(\bx)$.

In geometrically linear  dynamical systems, the convex function $T(\by)$ is the kinetic energy
and  $\UU(\bx) $  represents the total potential of the system.
 In this case,
the Lagrangian form $\Pi(\bx)= \TT(\DD \bx) - \UU(\bx)$   is the so-called {\em total action},
 which is a d.c. (difference of convex)  function.
 Since the  mixed  Lagrangian form $L(\bx, \by^*)$ is no longer a saddle function,
  the well-known Hamiltonian
\[
H(\bx, \by^*) = \la \DD \bx ; \by^* \ra - L(\bx, \by^*) = \TT^*(\by^*) + \UU(\bx)
\]
was introduced, which  is convex and  has been  extensively used in dynamical
  systems.
The Euler-Lagrange equations (\ref{eq-el}) is  equivalent to the
  well-known
 {\em canonical Hamiltonian equations}
  \eb
  \DD \bx = \nabla_{\by^*} H(\bx, \by^*) , \;\; \DD^* \by^* = \nabla_{\bx}  H(\bx, \by^*) .
  \ee
Actually, although the Lagrangian is not a saddle function in convex Hamiltonian systems,
 it is a so-called {\em super-critical function} \cite{GaoBook}, and if
 the total potential $\UU(\bx)$  is strictly convex on $\calX_a \subset \real^n$ such that
 its Legendre conjugate $\UU^*(\bx^*)$ can be uniquely defined on $\calX^*_a$,
 then
 the canonical dual action of $\Pi(\bx)$ can still be defined by
 \[
 \Pi^*(\by^*) = \max \{  L(\bx, \by^*) | \; \; \bx \in \calX_a \} = \UU^*(\DD^* \by^*) - \TT^*(\by^*)
 \]
 on $\calY^*_s = \{ \by^* \in \calY^*_a | \; \DD^* \by^* \in \calX^*_a \}$, which is also a d.c. function.
Therefore, instead of mono-duality in static systems,
 the convex Hamiltonian system is controlled by the so-called {\em bi-duality theory}.
  \begin{theorem}[Bi-Duality Theorem]
If  $(\barbx, \barby^*)$ is a critical point of the Lagrangian
  $L(\bx, \by^*)$, then $\barbx$ is a critical point of $\Pi(\bx)$, $\barby^*$ is a critical point of $\Pi^*(\by^*)$
  and   $ \Pi(\barbx) = L(\barbx, \barby^*) = \Pi^*(\barby^*)$.
  Moreover, if $n=m$,  we have   either
  \eb
 \Pi(\barbx) =  \max_{\bx \in \calX_k} \Pi(\bx) \;\;  \Leftrightarrow \;\;   \max_{\by^* \in \calY^*_s} \Pi^*(\by^*) = \Pi^*(\barby^*)
 \ee
 or
 \eb
 \Pi(\barbx) =  \min_{\bx \in \calX_k} \Pi(\bx)  \;\;  \Leftrightarrow \;\;   \min_{\by^* \in \calY^*_s} \Pi^*(\by^*) = \Pi^*(\barby^*).
 \ee
 \end{theorem}

This bi-duality is actually a special case of the triality theory in
geometrically linear systems, which was originally presented in Chapter 2 \cite{GaoBook}
for one-dimensional dynamical systems with a simple proof.
This bi-duality reveals a stable periodical property in convex Hamiltonian systems.

\subsection{Geometrically nonlinear systems and canonical duality}
Problems in geometrically nonlinear systems  are usually nonconvex.
Due to the fact that the geometrically linear operator $\DD:\real^n \rightarrow  \real^m$
 can not change the convexity of the objective function,
  if  $W(\bx)$ is nonconvex and  $W(\bx) = \TT(\DD \bx)$,
  the function $\TT(\by)$
 is still nonconvex and  its Legendre conjugate $T^*(\by^*)$ can not be uniquely defined\cite{sewell}.
 It turns out that    traditional Lagrangian duality  theory can not be applied  directly in this case.
 Although the Fenchel  conjugate
 $T^\sharp (\by^*) = \sup\{ \la \by ; \by^* \ra - T(\by) | \; \by \in \calY_a\}$ can be uniquely defined,
 the function
 \eb
{\L}(\bx, \by^*) = \la D \bx ; \by^* \ra - T^\sharp(\by^*) - U(\bx)
 \ee
 is not the traditional Lagrangian form and the associate
   saddle min-max duality theory will
 produce the so-called duality gap in nonconvex optimization.

Actually, in terms of $\UU(\bx) = \la \bx, \bff \ra - \half \la \bx, A \bx \ra$, the
total complementary function $\Xi(\bx, \bvsig) $ defined by (\ref{Complementary}) can be written as
\eb
\Xi(\bx, \bvsig) = \la \Lam(\bx) ; \bvsig \ra - \VV^*(\bvsig) - \UU(\bx).
\ee
Comparing  this $\Xi(\bx, \bvsig)$
with either ${\L}(\bx, \by^*)$  or the mixed Lagrangian form $L(\bx, \by^*)$
we can see that
the fundamental difference between the canonical duality  theory and  other methods 
is  the canonical transformation
$W(\bx) = \VV(\Lam(\bx))$ instead of
the  linear transformation $W(\bx) = \TT(\DD \bx)$ used in many other  duality theories,
 including the Fenchel-Moreau-Rockafellar duality.
In real applications,  if the quadratic function $\UU(\bx)$
is  nonconvex,  the mixed  Lagrangian form  $L(\bx, \by^*)$ is  nonconvex in $\bx$ since $\DD$ is linear.
 However,   the total complementary function
$\Xi(*, \bvsig):\calX \subset \real^n \rightarrow \real$ is always
  convex for  $\bvsig \in \calS^+_a$ and concave for  $\bvsig \in \calS^-_a$
  due to the geometrically nonlinear operator $\Lam(\bx)$ and its canonical dual variable $\bvsig$.
  Therefore, $\Xi(\bx,\bvsig)$ was also called the {\em nonlinear Lagrangian} in \cite{GaoBook}
  and the {\em extended Lagrangian} in \cite{gao-opt03}. If the geometrical operator $\Lam(\bx)$ is quadratic
  and objective, the so-called $\Lam$-transformation \cite{gao-opt03}
  \eb
  \UU^\Lam(\bvsig) = \sta \{ \la \Lam(\bx) ; \bvsig \ra - \UU(\bx) | \;\; \bx \in \calX \}
  \ee
  is actually the {\em pure complementary gap function} which is obtained from the complementary gap function
  $G_{ap}(\bx, \bvsig) = \half \la  \bx , G(\bvsig) \bx \ra $ by using the analytical solution form
  $\bx = [G(\bvsig)]^{-1} F(\bvsig)$.

The geometrical nonlinearity in continuum physics means large deformation (far from equilibrium states),
which usually leads to bifurcation in static systems \cite{santo-gao} and
chaos in dynamical systems \cite{Gao-amma03}.
Therefore,  geometrically nonlinear systems are usually nonconvex.
This is the reason why the  geometrical nonlinearity was emphasized in
the title of Gao and Strang's original work \cite{GaoStrang89}, although the system they studied is
fully nonlinear and governed by a nonconvex/nonsmooth total (super) potential functional
\eb
  \Pi(\uu) = \barW(\barLam (\uu) ) + \barF(\uu)  ,
  \ee
  where $\barW(e)$ is called the  stored   energy, which is a  canonical function(al)
  such that the  constitutive law $e^* = \partial \barW(e)$ is invertible on its effective domain;
   while $\barF(u)$ is an external energy,
  which must be linear on the statically admissible space
   such that its \G   derivative $\partial \barF(u) = - \baru^*$ leads to the
  external force (source) field  (under the sign convention).
  The  geometrically nonlinear operator $e = \barLam(\uu)$
  in Gao and Strang's work should be an objective   measure
  in order to satisfy certain well-known  deformation laws (see Chapter 6, \cite{GaoBook}).
Therefore, the complementary  gap function $G_{ap}(u, e^*)$ was naturally introduced.
This objective function  lays a foundation for the triality theory.

Oppositely, in a  recent paper entitled
 ``Some remarks concerning Gao-Strang's complementary gap function"
  by Voisei and Z{a}linescu \cite{vz1}, they
choose quadratic functions as the external energy $\barF(u)$ (see Examples 2, 4 and 5 in \cite{vz1}),
 and   piecewise linear function
 (see Example 1 in \cite{vz1}) as the stored energy,
 they concluded: ``About the (complementary) gap function one can conclude that it is useless at least in the current context".
 Clearly, the piecewise linear function is not objective and cannot store energy;
 while for those quadratic functions  $\barF(u)$ they listed, the dual variable  $u^* = \partial  \barF(u)$
 depend on the configuration $u$. Such force field is called
  {\em follower force}. In this case, the system is not conservative and traditional variational methods do not apply.
 Unfortunately, similar counter-examples and conclusions are repeatedly presented  in many other papers (see \cite{vz-dcds} and references cited therein).

Actually, in order to  study
nonconvex variational problems in dissipative systems subjected to follower force field,
 a so-called rate variational  method and the associated dual extremum principle were proposed in 1990 \cite{gao-onat}.
 Also, Gao and Strang's work has been extended to general nonconvex dynamical systems
 to allow $\barF(u)$ as a quadratic function, but  notations  were  changed  (see  \cite{gao-opt03,Gao-amma03}).
In fact,
 if we let  $\barLam(\bx) = \{ \Lam(\bx),
 \half \la \bx, A \bx \ra \}$ and $\barW(\barLam(\bx)) = \VV(\Lam(\bx)) + \half \la \bx, A \bx \ra$,
  the general nonconvex problem $(\calP)$ studied in this paper
 is simply a finite dimensional version of the
  Gao and Strang's  general work in large deformation theory.
  This method has been repeatedly used in many Gao's papers (see \cite{gao-ogden,yau-gao}).
Particularly, if $\barLam (u)$ is a Cauchy-Riemann strain measure, then
  \eb
  \Xi(u, e^*) = \la \barLam(u) ;  e^* \ra - \barW^*(e^*) + \barF(u)
  \ee
is the well-known {\em Hellinger-Reissner complementary energy} in finite deformation theory\footnote{The Hellinger-Reissner energy
 was first proposed by Hellinger in 1914. After the external energy $\barF(u)$ and the
boundary conditions in the statically admissible
 space $\calU_k = \{ u\in \calU_a | e = \barLam(u) \in \calE_a \}$
 were fixed by Reissner in 1953,
the associated variational statement has been  known as the Hellinger-Reissner principle.
However, the extremality condition of this  principle was an open problem, and also
  the existence of pure complementary variational principles  has been a
well-known debate existing for over several decades in large deformation mechanics (see \cite{li-gupta}).
This open problem was partially solved by  Gao and Strang's work and completely solve by
the triality theory. While the
pure complementary energy principle was  formulated by Gao in 1999 \cite{gao-mecc}.}.
Furthermore, if the complementary energy $\barW^*(e^*)$ is replaced by
 $\la e ; e^* \ra - \barW(e )$, the total complementary energy $\Xi(u, e^*)$ can be written in
 the so-called {\em pseudo-Lagrangian} (it was denoted as $L_p(u, e^*, e)$ in \cite{GaoStrang89})
 \eb
\Xi_{hw}(u, e^*, e) = \barW(e) + \la \barLam(u) - e ; e^* \ra + \barF(u),
\ee
and we have
\[
\Xi(u, e^*) = \sta \{ \Xi_{hw}(u, e^*, e) | \; e\in \calE_a \}.
\]
 In large deformation mechanics, $\Xi_{hw}(u, e^*, e)$ is  called the
 Hu-Washizu generalized potential energy, proposed independently by
 Hai-Chang Hu in 1954 and K. Washizu in 1955.
 The  associated variational statement is
 the  well-known {\em Hu-Washizu principle}, which has important applications in
 computational   mechanics of thin-walled structures, where the geometrical equation
 $e = \barLam(u)$ is usually proposed by certain geometrical hypothesis
\cite{gao-cheung,Gao-Sherali-AMMA09}.

It has been  emphasized in many papers that
 the key step in the canonical duality theory is to choose
a  geometrically reasonable measure $\bxi = \Lam(\bx)$.
It was shown in \cite{gao-yang} that for a given nonconvex variational problem,
the choice of  $\Lam(\bx)$ may not be unique and different geometrically admissible operators
could lead to different canonical dual problems.
But all   these  canonical dual problems must be equivalent in the sense that they have the
same set of solutions.
Also for complex systems, two type of sequential canonical transformations were proposed
(see Chapter 4, \cite{GaoBook}).
 By the fact that the objectivity and canonical duality are fundamental to all natural systems,
    for any given real  problem, as long as the  geometrical operator $\Lam(\bx)$
 can be chosen correctly such that
 the nonconvex objective function(al) can be recast by adopting a canonical form  $
 W(\bx) = V(\Lam(\bx))$, the
 canonical duality theory can be used to establish elegant  theoretical results and to develop
efficient  algorithms for robust computations.
   The  triality
theory reveals an intrinsic duality pattern in
 nonconvex systems and
 should play important roles not only for
 solving a large class of challenging problems in
 nonconvex analysis and global optimization, but also for understanding, modeling,   and simulation of complex  systems.

\section{Conclusion Remarks}
Motivated by an open problem on the double-min duality in the triality theory that was left
unaddressed since 2003, we have presented a mathematically rigorous proof for this theory based
on the elementary linear algebra.
 Our results show that the triality theory   holds strongly in
the tri-duality form if the primal and its dual problems have the
  same dimension. Otherwise, both the canonical  min-max and
double-max duality statements  hold strongly, but
 double-min duality statement   holds weakly in a super-symmetric
form. Additionally, a weak saddle-duality theory is proposed, which shows that
when the complementary gap function $G_{ap}(\bx, \bvsig)$
 is negative, either the primal problem $(\calP)$ (only if $m< n$)
or its canonical dual $(\calP^d)$ (only if $m > n$) could have saddle critical solutions.
Therefore, this seven years old open
problem is now solved completely
and the triality theory is presented in an elegant form as expected.

The method adopted in this paper can be generalized for more general constrained global optimization problems.
As it is mentioned in Remark 3 that the
primal problem $(\calP)$ could be NP-hard if its canonical dual has no critical point in  $\calS^+_a$.
Also, the extremality conditions for those critical points $(\barbx, \barbvsig)$
are still unknown if the Hessian matrix $G(\barbvsig)$ of the gap function is indefinite.
Although a general theorem on the existence and uniqueness  of the canonical dual solution in $\calS^+_a$
was proposed in \cite{gao-cace09}, and some perturbation methods
were discussed in \cite{gao-ruan-jogo10},
detailed quantitative study on these  topics
is fundamentally important and critical for understanding  and solving NP-hard problems.\\


\noindent {\bf Acknowledgements}
The authors are gratefully indebted with  Professor Hanif Sherali at Virginia Tech for
his detailed remarks and important suggestions.
This paper has benefited from three anonymous referees' constructive comments.
The main results of this paper were announced at the 2nd World Congress  of Global Optimization,
July 3-7, 2011, Chania, Greece.
 David Gao's research is supported by US Air Force
Office of Scientific Research under the grant AFOSR FA9550-10-1-0487.
 Changzhi Wu was supported by National Natural
Science Foundation of China under the grant \# 11001288, the Key
Project of Chinese Ministry of Education under the  grant \#
210179, SRF for ROCS, SEM, Natural Science Foundation Project of
CQ CSTC under the  grant \# 2009BB3057 and CMEC under the grant \#
KJ090802.


\begin{thebibliography}{99}

\bibitem{arnold} Arnold, V.I.: On teaching mathematics,
{\em Uspekhi Mat. Nauk.} 53 (1998), no. 1, 229—234;
English translation: {\em Russian Math. Surveys} 53 (1998), no. 1, 229—236.

\bibitem{fang-gaoetal07} Fang, S.C., Gao D.Y., Sheu R.L. and Wu S.Y.:  Canonical
dual approach for solving 0-1 quadratic programming problems, {\em
J. Ind. and Manag. Optim.} 4, 125-142 (2008).

\bibitem{Matrix} Filly, J.A. and Donniell, E.F.: The Moore-Penrose generalized
inverse for sums of matrices, SIAM J. Matrix. Anal. Appl., Vol.21,
No. 2, pp. 629-635, (1999).

\bibitem{Floudas1}Floudas C.A.: Deterministic Global Optimization: theory,
methods and applications. Kluwer Academic, Dordrecht (2000).

\bibitem{Floudas2} Floudas, C.A.: Systems approaches in bioinformatics and
computational genomics, Challenges for the Chemical Sciences in
the 21th Century, Information and Communications Workshop,
National Research Council of the National Academies, National
Academies Press, pp 116¨C125 (2003).



\bibitem{gao-tri96}Gao, D.Y.:
 Post-buckling analysis and anomalous dual variational problems in nonlinear beam theory.
  {\em Applied Mechanics in Americans, Proc. of the Fifth Pan American
  Congress of Applied Mechanics}, L.A. Godoy, L.E. Suarez (Eds.), Vol. 4.
   The University of Iowa, Iowa city, August (1996).

 \bibitem{gao-amr} Gao, D.Y.:
    Dual extremum principles in finite deformation theory
    with applications to  post-buckling analysis of extended nonlinear beam theory.
    {\em Applied Mechanics Reviews},  50 (11),    S64-S71 (1997).

\bibitem{gao-mecc}   Gao, D.Y. (1999).    General Analytic Solutions and
  Complementary Variational Principles for Large Deformation
Nonsmooth Mechanics. \emph{Meccanica}  {34},
169-198.

\bibitem{GaoBook} Gao, D.Y.:
{\em  Duality Principles in Nonconvex Systems: Theory, Methods and
Applications.} Kluwer Academic, Dordrecht, 2000.

\bibitem{GaoJOG00} Gao, D.Y.: Canonical dual transformation method and
generalized triality theory in nonsmooth global optimization. J.
Glob. Optim. 17(1/4), 127--160 (2000).



\bibitem{gao-opt03} Gao, D.Y.: Perfect duality theory and complete solutions
to a class of global optimization problems. {\em Optim.} 52(4--5),
467--493(2003)

\bibitem{Gao-amma03} Gao, D.Y.:  Nonconvex semi-linear problems and canonical dual solutions.
In: Gao, D.Y., Ogden, R.W.
(eds.) {\em Advances in Mechanics and Mathematics,} vol. II, pp. 261–312. Kluwer Academic, Dordrecht
(2003)


\bibitem {gao-jimo07} Gao, D.Y. (2007).
Solutions and optimality  to box constrained nonconvex
minimization problems {\em J. Indust.  and Manage.
Optim.}, 3(2), 293-304.



\bibitem{gao-cace09} Gao, D.Y.: Canonical duality theory: theory, method, and
applications in global optimization. {\em Comput. Chem.} 33,
1964-1972, (2009).

\bibitem{gao-cheung}
Gao, D.Y. and Cheung, Y.K.: On the extremum complementary energy principles for nonlinear elastic shells,
{\em Int. J. Solids \& Struct.,} 26 (1989), pp. 683-693.



\bibitem{gao-ogden} Gao, D.Y. and Ogden, R.W.: Multiple solutions to non-convex
variational problems with implications for phase transitions and numerical
computation. {\em Q. J. Mech. Appl. Math,} Vol. 61. No. 4, 497-522.


\bibitem{gao-onat}
Gao, Y. and Onat, E.T.: Rate variational extremum principles for finite elasto- plasticity.
{\em Appl. Math. Mech.,} 11 (1990), 7, pp. 659-667.

\bibitem{gao-ruan-jogo10} Gao, D.Y. and Ruan, N.:  Solutions to quadratic minimization
problems with box and integer constraints. {\em J. Glob. Optim.}
47(3): 463-484, (2010).





\bibitem{gao-ruan-pardalos} Gao, D.Y., Ruan, N, and Pardalos, P.M. (2010).
Canonical dual solutions to sum of fourth-order
polynomials minimization problems with applications to
sensor network localization, in {\em Sensors: Theory, Algorithms and Applications},
P.M. Pardalos, Y.Y. Ye, V. Boginski, and C. Commander (eds). Springer.

\bibitem{Gao-Sherali-AMMA09} Gao, D.Y. and Sherali, H.D.:  Canonical duality:
Connection between nonconvex mechanics and global optimization, in
{\em Advances in Appl.  Mathematics and Global Optimization}, 249-316,
Springer (2009).


\bibitem{GaoStrang89} Gao, D.Y., Strang, G.: Geometric nonlinearity:
Potential energy, complementary energy, and the gap function.
Quart. Appl. Math. 47(3), 487--504 (1989).

\bibitem{gao-watsonetal}Gao, D.Y.,   Watson, L.T., Easterling,  D.R., Thacker, W.I.
and  Billups, S.C.:
 {Solving the canonical dual of box- and integer-constrained nonconvex quadratic programs via a deterministic direct search algorithm},
{\em Optimization Methods \& Software, } (2011)
 DOI: 10.1080/10556788.2011.641125

\bibitem{GaoWu} Gao, D.Y. and Wu, C.Z.: On the triality theory in global optimization.
 Announced at http://arxiv.org/abs/1104.2970,  arXiv:1104.2970v1  (2010).

\bibitem{gao-wu-jimo11}Gao, D.Y. and Wu, C.Z.:
On the triality theory for a quartic polynomial optimization problem,
 {\em J. Industrial and Management Optimization}, 8(1):229-242, (2012).


\bibitem{gao-yang}  Gao, D.Y. and Yang, W.-H.:
Multi-duality in minimal surface type problems,
{\em Studies in Appl. Math.}, 95,   127-146 (1995).

\bibitem{Gao-Yu2008}  Gao, D.Y. and Yu, H.F.:
  Multi-scale modelling and canonical dual finite element method in
  phase transitions of solids.
  {\em Int.  J.  Solids and Structures}, 45:3660--3673 (2008).

\bibitem{holz}Holzapfel, G.A.:
 {\em Nonlinear Solid Mechanics: A Continuum Approach for Engineering.} Wiley, (2000).
 ISBN 978-0471823193.

\bibitem{Ogden}Holzapfel, G.A. and Ogden R.W.(eds):
{\em Mechanics of biological tissue,}  vol 12. Springer, Heidelberg, 522 pp.(2006)

\bibitem{landau-lif} Landau, L.D. and Lifshitz, E.M. (1976). {\em Mechanics}. Vol. 1 (3rd ed.). Butterworth-Heinemann. ISBN 978-0-750-62896-9.

\bibitem{li-gupta}
Li, S.F. and Gupta, A.: On dual configuration forces, {\em J. of Elasticity, }
84:13-31 (2006).

\bibitem {moreau} Moreau, J.J.:
 La notion de sur-potentiel et les liaisons unilat\'{e}rales en
\'{e}lastostatique, {\em C.R. Acad. Sc. Paris,} 267 A, (1968), 954-957.

\bibitem{ogden} Ogden, R.W.: {\em Non-Linear Elastic Deformations}, Dover Publications, 1 edition (July 7 1997), 544pp.

\bibitem{r-g-j}Ruan, N.,  Gao,  D.Y.,  and Jiao, Y.:
Canonical dual least square method for solving general nonlinear
systems of quadratic equations, {\em Computational Optimization and Applications},
Vol 47, 335-347 (2010).

\bibitem{santo-gao}Santos, H.A.F.A. and Gao D.Y.:  Canonical dual finite element method for solving post-buckling problems of a large deformation elastic beam. {\em
     Int. J. Nonlinear Mechanics,} 2011,  doi:10.1016/j.ijnonlinmec.2011.05.012
\bibitem{sewell}
Sewell, M.J.: {\em Maximum and Minimum Principles: A unified approach, with applications}.
 Cambridge University Press, 1987, 468pp.

\bibitem{silva-gao-jmaa1} Silva, D.M.M and Gao, D.Y.:
Complete solutions and triality theory to a nonconvex optimization problem with double-well potential
in $\real^n$,  2011, 	arXiv:1110.0285v1



\bibitem{gs86}  Strang, G.:  {\em Introduction to Applied
Mathematics}, Wellesley-Cambridge Press, 1986, 758 pp.


\bibitem{vz-dcds}  Strugariu, R.,  Voisei, M.D., and Zalinescu, C.: Counter-examples in bi-duality, triality and tri-duality.
    {\em Discrete and Continuous Dynamical Systems - Series A},
     Volume 31, Number 4, December 2011.


      \bibitem {tonti72} Tonti, E.:
On the mathematical structure of a large
class of physical theories,
{\em   Accad. Naz. dei  Lincei},  Serie VIII,
{ LII},  49-56 (1972).

\bibitem{vz1} M.D. Voisei, C. Zalinescu: Some remarks concerning Gao-Strang's
 complementary gap function, {\em Applicable Analysis,}
 DOI: 10.1080/00036811.2010.483427, 2010.

\bibitem{vz-jogo}M.D. Voisei, C. Zalinescu:
Counterexamples to some triality and tri-duality results.
{\em J. Global Optimization},  49:173–183, 2011.

\bibitem{wang-etal} Wang, Z.B., Fang, S.C., Gao, D.Y., and Xing, W.X.:
Canonical dual approach to solving the maximum cut problem, to
appear in {\em J. Glob. Optim.}


 \bibitem{yau-gao} Yau, S.-T. and Gao,  D.Y.:
 Obstacle problem for von Karman equations, {\em Adv. Appl. Math.}, 13,
123-141 (1992).

\bibitem{zhang-gao}Zhang J., Gao, D.Y. and Yearwood, J.:
 A novel canonical dual computational approach for prion AGAAAAGA amyloid fibril molecular modeling. {\em Journal of Theoretical Biology}, 284,  149-157  (2011). doi:10.1016/j.jtbi.2011.06.024
\end{thebibliography}
\end{document}